\documentstyle[12pt,amssymb]{article}
\newcommand{\eq}{\begin{equation}}
\newcommand{\en}{\end{equation}}

\newtheorem{theorem}{Theorem}
\newtheorem{lemma}{Lemma}
\newtheorem{corollary}{Corollary}

\newtheorem{example}{Example}
\newtheorem{remark}{Remark}

\def\BP{{\Bbb P}}
\def\BP{{\Bbb P}}
\def\BR{{\Bbb R}}

\def\BE{{\Bbb E}}
\def\BZ{{\Bbb Z}}
\def\BC{{\Bbb C}}

\def\cs{c^\ast}
\def\Bs{B^\ast}

\def\Ms{M^\ast}
\def\phis{\phi^\ast}
\def\gs{g^\ast}
\def\Xs{X^\ast}
\def\sigmas{\sigma^\ast}
\def\Xs{\tilde{X}}
\def\Ys{\tilde{Y}}
\def\sigmas{\tilde{\sigma}}

\def\calF{{\cal F}}

\newcommand{\ignore}[1]{}

\begin{document}

\title{Asymptotic enumeration and logical limit laws
for expansive multisets and selections}
\author{Boris L. Granovsky\thanks{
Department of Mathematics, Technion-Israel Institute of
Technology, Haifa, 32000, Israel,
e-mail:mar18aa@techunix.technion.ac.il
 }
 \ and  Dudley
Stark\thanks{School of Mathematical Sciences, Queen Mary,
University of London, London E1 4NS, United Kingdom,
e-mail:D.S.Stark@maths.qmul.ac.uk\newline\newline 2000 {\it
Mathematics Subject Classification}: 6 0C05 (primary),
05A16 (secondary)} }

\date{May 4, 2005}

\maketitle

\begin{abstract}
Given a sequence of integers $a_j, \ j\ge 1,$
a multiset is a combinatorial object
composed of unordered
components, such that there are exactly $a_j$
 one-component multisets  of size $j.$
  When $a_j\asymp j^{r-1} y^j$ for some $r>0$,
$y\geq 1$, then  the multiset is called {\em expansive}. Let
$c_n$ be the number of multisets  of total size  $n$. Using
a probabilistic approach, we prove  for expansive multisets
that $c_n/c_{n+1}\to 1$ and that $c_n/c_{n+1}<1$ for large
enough $n$. This allows us to prove Monadic Second Order
Limit Laws for expansive multisets. The above results are
extended to a class of expansive multisets with
oscillation.

Moreover, under  the condition $a_j=Kj^{r-1}y^j + O(y^{\nu j}),$
where $K>0$, $r>0$, $y>1$, $\nu\in (0,1)$, we find an explicit
asymptotic formula for $c_n$. In a similar way we study the
asymptotic behavior of selections which are defined as
combinatorial objects composed of unordered components of
distinct sizes.

\end{abstract}
\vfill

\section{Summary and  Historical remarks}
  Given a sequence of integers $a_j\ge 0, \  j\ge 1$, a
multiset is a
 combinatorial object of finite total size composed of unordered
indecomposable components such that there are exactly $a_j$
single component multisets  of size $j.$ There is no
restriction on the number of times a component may appear
in the multiset.In view of this, let

$$ \Omega_n=\left\{
\vec{\eta}=(\eta_1,\eta_2,\ldots,\eta_n):
\sum_{j=1}^nj\eta_j=n {\rm \ and \ }\eta_j\geq 0 {\rm \ for
\ all \ }j\right\} $$ be the set of unordered integer
partitions of an integer $n$. Any multiset of total size
$n$ has a component count vector $
\vec{\eta}=(\eta_1,\eta_2,\ldots,\eta_n)$ contained in
$\Omega_n$. Here $\eta_j=\eta_j^{(n)}$ is the number of
components of size $j$ in the multiset of size $n$
considered.(For more details about multisets see
e.g.\cite{ABT,Bur}).

 Let $c_n$ be the number of
multisets of size $n$ determined by the above parameters
$a_j, \ j=1,\ldots, n$. We will prove an analytic identity
which will be used to extract information about the growth
of $c_n$, as $n\to \infty$. It follows from the
 definition of a multiset that the number of multisets
with a given component count vector $
\vec{\eta}=(\eta_1,\eta_2,\ldots,\eta_n)\in \Omega_n$
 is $$\prod_{j=1}^n{{a_j+\eta_j-1}\choose \eta_j}.$$
Hence, the number of multisets of size $n$ has the Euler
type generating function $g$:
\begin{eqnarray}\label{ogf}
g(x)&=&\sum_{n=0}^\infty c_n x^n\nonumber\\ & =&
1+\sum_{n=1}^\infty\sum_{\vec{\eta}\in\Omega_n}
\prod_{j=1}^n{{a_j+\eta_j-1}\choose \eta_j}x^{j\eta_j}\nonumber\\
&=& \prod_{j=1}^\infty\left(1-x^j\right)^{-a_j}, \quad
\vert x\vert <1.
\end{eqnarray}
We observe that combinatorial objects that are given by the
generating function (\ref{ogf}) are also called weighted
partitions(\cite{Br}).

 The truncated generating function
\eq\label{trunc1} g_n(x)=\prod_{j=1}^n\left(1-x^j\right)^{-a_j},
\quad  \vert x\vert<1 \en has Taylor expansion
$g_n(x)=\sum_{k=0}^\infty c_{k,n} x^k$, where $c_{k,n}=c_n$ for
$k\leq n$. For a fixed $n$, the series expansion of the function
$g_n(x)$ converges for all $x\in\BC$, $|x|<1$. We set $$
x=x(\sigma,\alpha)=e^{-\sigma+2\pi i \alpha} $$ for some real
numbers $\alpha$ and $\sigma$. Then we have \eq\label{fourier}
\int_0^1g_n(x)e^{-2\pi i\alpha n}d\alpha=
\int_0^1\left(\sum_{k=0}^\infty c_{k,n}e^{-k\sigma + 2\pi
i\alpha(k-n)}\right) d\alpha =c_ne^{-n\sigma}, \en where we have
used $$ \int_0^1e^{2\pi i \alpha m}d\alpha= \left\{
\begin{array}{l l}
1,& {\rm if \ }m=0,\\
0,& {\rm if \ }m\neq 0, \ \ m\in\BZ.
\end{array}
\right.
$$
Substituting (\ref{trunc1}) into (\ref{fourier}) gives the
desired identity  w.r.t.~the free parameter $\sigma\in \BR
$ :

 \eq\label{ident1} c_n=e^{n\sigma} \int_0^1
\prod_{j=1}^n \left( {1-e^{-j\sigma}e^{2\pi i\alpha
j}}\right)^{-a_j} e^{-2\pi i\alpha n} d\alpha. \en
 A more
probabilistic interpretation of (\ref{ident1}) can be given. We
have
\begin{eqnarray}\label{rewrite1}
c_n&=&e^{n\sigma} \prod_{j=1}^n \left(
{1-e^{-j\sigma}}\right)^{-a_j} \int_0^1 \prod_{j=1}^n \left(
\frac{1-e^{-j\sigma}}{1-e^{-j\sigma}e^{2\pi i\alpha
j}}\right)^{a_j} e^{-2\pi i\alpha n} d\alpha\nonumber\\
&=&e^{n\sigma} \prod_{j=1}^n \left(
{1-e^{-j\sigma}}\right)^{-a_j} \int_0^1 \phi(\alpha)\, e^{-2\pi
i\alpha n} d\alpha,
\end{eqnarray}
where $\phi(\alpha)$ is given by \eq\label{phidef}
\phi(\alpha)=\prod_{j=1}^n \phi_j(\alpha),\quad \alpha\in\BR, \en
for functions $\phi_j$ defined by \eq\label{phikdef}
\phi_j(\alpha)= \left(\frac{1-e^{-j\sigma}}{1-e^{-j\sigma}e^{2\pi
i\alpha j}}\right)^{a_j},\quad\alpha\in\BR. \en
 Using the combinatorial identity in Lemma 2.46 in \cite{Bur}
it is easy to see that for any $\sigma>0$, $\phi_j$ is the
characteristic function of a random variable $X_j$  given by
\eq\label{xyz} \BP(X_j=jl)={{a_j+l-1}\choose
l}\left(1-e^{-j\sigma}\right)^{a_j} e^{-lj\sigma}, \ \
l=0,1,2,\ldots. \en Consequently,
 $\phi(\alpha)$ is the characteristic function of $Y_n:=\sum_{j=1}^n
X_j$, where the $X_j,\ j=1,\ldots,n$ are assumed to be
independent. Therefore, \eq\label{zz}
\int_0^1\phi(\alpha)e^{-2\pi i\alpha n}d\alpha
=\BP\left(Y_n=n\right).\en

Combining (\ref{zz}) with (\ref{rewrite1}) we arrive  at
the desired representation of $c_n,$ which is in the core
of the probabilistic method  suggested by Khintchine in
1950-s ( \cite{K}, Chs IV, V) for  asymptotic enumeration
in the urn models of statistical mechanics.The history
related to the method is discussed in \cite{FG2}.  We note
that Khintchine-type representations were subsequently
rediscovered in  independent ways by many authors in a
variety of seemingly unrelated contexts. In particular,
observe that (\ref{rewrite1}) can be derived from equation
(134) of \cite{AT}, the latter being based  on the
conditioning relation (see \cite{ABT}). In conclusion, the
recent paper \cite{MK} should be mentioned which treats
probabilistic methods for enumeration as transforms of
generic random variables (in our setting $Z_j$) into
specially constructed independent random variables (in our
setting $X_j$).

It follows from (\ref{xyz}) that the r.v. $j^{-1}X_j$ is
negative binomially distributed with parameters $a_j$ and
$e^{-\sigma j}, \ \sigma>0$. This  produces the formula for
the expectation of the sum $Y_n$: \eq\label{Mdef} M_n:=\BE
Y_n= \sum_{j=1}^n\frac{ja_je^{-j\sigma}}{1-e^{-j\sigma}}.
\en Further on, except for Lemma \ref{phiexp} in Section 2,
we will assume that $ \sigma=\sigma_n>0$ is chosen in the
unique way so that \eq \label{12} M_n=n. \en The fact that
$\sigma$ can be chosen in such a way follows from observing
that $M_n$ decreases monotonically from $\infty$ to $0$ as
$\sigma$ ranges from $0$ to $\infty$, while $n$ is
fixed. The variance of $Y_n$ is \eq\label{Bdef} B_n^2:={\rm
Var}(Y_n)=
\sum_{j=1}^n\frac{j^2a_je^{-j\sigma_n}}{(1-e^{-j\sigma_n})^2}.
\en

 We will study the probability in (\ref{zz})
analytically and thereby obtain information about
asymptotic behavior of $c_n, $ as $n\to \infty$.

It is natural to suppose that, under some asymptotic
conditions on the parameters $a_j, \ j\ge 1$, a local limit
law should hold for $Y_n$ in (\ref{zz}). Asymptotic
enumeration of multisets using this approach was apparently
initiated in \cite{FP}, which was preceded by works of
Freiman ( see \cite{P}) on the development of Khintchine's
method. Note that asymptotics  of different statistics
related to integer  partitions (the case $a_j=1,\ j\ge 1$)
was studied by many authors(see e.g.\cite{FR,M}).

 In this paper we will initially
assume that \eq\label{exp} a_j\asymp j^{r-1}y^j,\quad j\to \infty,
\quad r>0,\quad y\ge 1, \en where we write $a_j\asymp b_j$ for
sequences $a_j$, $b_j$ when there exist constants $D_1,D_2>0$ such
that $D_1\leq a_j/b_j \leq D_2$ for all $j\ge 1$. Although for
$y>1$, the rate of growth of $a_j$ is exponential (but not
polynomial) such multisets will be called, following \cite{BG},
{\em expansive}. This is in view of Bell-Burris Lemma 5.2 in
\cite{BB} which tells us that for $y>1$, the asymptotic behavior
of the number of multisets with $a_j\sim j^{r-1}y^j$ is the
same as that of assemblies with $a_j\sim j^{r-1}$. (Here and in
what follows $a_n\sim b_n$ for sequences $a_n$, $b_n$ means that
$\lim_{n\to\infty}a_n/b_n=1$).

Provided the parameters  $a_j$ satisfy (\ref{exp}), we will prove
that the normal limiting law for the sum $Y_n$ holds, leading
to $\BP(Y_n=n)\sim (2\pi B_n^2)^{-1/2}$ in (\ref{zz}). This local
limit law, the definition of $\sigma_n$, and (\ref{rewrite1}),
will imply our Theorem \ref{main} below. A result implying Theorem
\ref{main} when in (\ref{exp}) $y=1$ and $\asymp$ is replaced by
$\sim$, was obtained by Richmond in \cite{R} and applied for
enumeration of partitions of $n$ into primes or powers of primes.
The first probabilistic proof of the Hardy-Ramanujan formula for
partitions ( the case $a_j=1, \ j\ge 1$) was given by Freiman in
1950-s (see Section 2.7 of \cite{P}.) Quite recently, a different
probabilistic proof of this formula was published in \cite{Ba}. A
comprehensive study of the asymptotics of integer partitions was
made in \cite{Pi}.

Theorem \ref{main} appears to be new for $y>1$.
Note that throughout the paper we assume, unless it is said
otherwise, that all asymptotic expressions are valid under
$n\to \infty$.

\begin{theorem}\label{main}
Assume that (\ref{exp}) holds. Then the number of multisets
is asymptotically \eq\label{mainasymp}
c_n\sim\frac{e^{n\sigma_n}}{\sqrt{2\pi B_n^2}} \,
\prod_{j=1}^n \left( {1-e^{-j\sigma_n}}\right)^{-a_j}, \en
where  $\sigma_n$ is given by (\ref{12}).

\end{theorem}

We now formulate  an extension  of Theorem \ref{main}
to a wider than (\ref{exp}) class of parameters $a_j$. Corollary 1
below is an analog of Corollary 1 in \cite{FG} for
expansive assemblies.

We write $q_1(n)\preceq\bullet(n)\preceq q_2(n)$,
 if there exist positive constants
$\gamma_1$, $\gamma_2$, such that $\gamma_1q_1(n)
\leq\bullet(n)\leq\gamma_2 q_2(n),
 n\ge 1 $ . For given $0<r_1\leq r_2$ and $y\geq 1$ define the set
$\calF(r_1,r_2,y)$ of parameter functions $a=a_j$, $j\ge 1$
satisfying the condition $$ j^{r_1-1}y^j\preceq a_j \preceq
j^{r_2-1}y^j, \ \ y\ge 1,\ j\ge 1. $$
\begin{corollary}\label{weaker1}
For an arbitrary $r>0$, $0<\epsilon\leq r/3$ and $y\geq 1$, the
conclusion of Theorem \ref{main} is valid for all parameter
functions $a\in\calF(2r/3+\epsilon,r,y)$.
\end{corollary}

It is interesting to  observe that in the case $y=1,$ our
condition $a\in\calF(2r/3+\epsilon,r,y)$ implies the condition
(i), p. 1084 of Richmond \cite{R}. This latter condition is
sufficient for the asymptotic formula for partitions of $n$
into primes developed in \cite {R}. Generally speaking, multisets
with $a\in\calF(r_1,r_2,y)$ may be called ``expansive with
oscillation''.

Theorem~\ref{main} and Corollary~\ref{weaker1} are proved in
Section 2. In Section 2 we also derive
 asymptotic estimates for $ \sigma_n$ and  $ B_n^2 $ that are used in
(\ref{mainasymp}).

A multiset satisfies a monadic second order  logical limit law if
the probability that a random representative of the multiset of
size $n$ satisfies a monadic second order sentence converges,
as $n\to \infty$. Compton \cite{Ca,Cb} showed that to prove that a
class of relational structures such as multisets satisfies a
monadic second order limit law, it suffices to know about the
growth of the number of structures $c_n$ of size $n$. The next
corollary from Compton's theorem was used in \cite{BB} to prove
logical limit laws.

\begin{theorem}\label{basic}{\rm\bf [Corollary 8.1 of \cite{BB}]}
Suppose that $a_j$ are the parameters of a multiset such that
$$
\frac{c_n}{c_{n+1}}\sim y^{-1} {\rm \ for \ some \ }y\geq 1.
$$
If $y>1$, then suppose further that there exists $N>0$ such that
$$
\frac{c_n}{c_{n+1}}\leq y^{-1} {\rm \ whenever \ }n>N.
$$
Then the multiset has a monadic second order logical limit law.
\end{theorem}

Based on  a Tauberian theorem  of Schur, Bell and Burris
(\cite{BB},Theorems 9.1 and 9.3) and Bell \cite{Bl} derived
general sufficient  conditions on the asymptotics of $a_j$
which imply the hypotheses of Theorem \ref{basic}. Note that
the condition on the $a_j$ obtained in \cite{BB} covers the
particular case $a_j\sim j^{r-1}y^j$, $y\geq 1$ of  (\ref{exp}).
Combining Theorem \ref{main}  with  the asymptotic estimates
in Section 2,  we prove in Theorem \ref{zeroone} below the
validity of the conditions of Theorem \ref{basic}, in the case
(\ref{exp}) that is not covered by the Bell-Burris sufficient
conditions.
\begin{theorem}\label{zeroone} Suppose that $a_j$ is a
sequence obeying the condition (\ref{exp}). Then the
corresponding multiset has a monadic second order logical limit
law.
\end{theorem}

A result similar to Theorem~\ref{zeroone} is obtained for logarithmic
structures in \cite{S}.

Moreover, we are able to weaken the condition
(\ref{exp}) of Theorem \ref{zeroone}:

\begin{corollary}\label{weaker2}
For all parameter functions $a\in\calF(2r/3+\epsilon,r,y),$
where $r>0,\ y\geq 1$ and $0<\epsilon\leq r/3,$

$$
\frac{c_n}{c_{n+1}}=y^{-1}\exp\left(-\delta_n+o\left(\delta_n\right)\right),
$$
where $\delta_n=\sigma_n-\log y\to 0$.
\end{corollary}

In particular, Corollary \ref{weaker2} implies that
$c_n/c_{n+1}\to y^{-1}, \ y\ge 1$, for
$a\in\calF(2r/3+\epsilon,r,1)$. A similar result was shown
in \cite{FG,FG2}  for certain reversible
coagulation-fragmentation processes. From an analytical
point of view  the latter processes are equivalent to
expansive assemblies (see \cite{FG2}).

In view of the above definition, we may consider multisets  as
unlabelled decomposable combinatorial structures. We call
labelled decomposable combinatorial objects assemblies, a term
used in \cite{AT}, see also \cite{ABT}. Sufficient conditions to
have monadic second order logical limit laws were given for both
multisets and assemblies in Theorem 6.6 of \cite{Ca}. Let $m_j$ be
the number of labelled components of size $j$ and let
$a_j=\frac{m_j}{j!}$. For assemblies the basic representation
(\ref{rewrite1}) becomes for an arbitrary choice of $\sigma$ $$
c_n=n!\,e^{n\sigma} \exp\left(\sum_{j=1}^n a_j e^{-j\sigma}\right)
\int_0^1 \psi(\alpha)\, e^{-2\pi i\alpha n} d\alpha, $$ where $$
\psi(\alpha)= \exp\left(\sum_{j=1}^na_je^{-j\sigma}(e^{2\pi\alpha
j}-1)\right),\quad r>0; $$ see (2.24) of \cite{FG} or (125) of
\cite{AT}. The method of proof of Theorem \ref{zeroone} and the
comment in the last paragraph gives monadic second order logical
limit laws for assemblies whenever $a\in\calF(2r/3+\epsilon,r,1)$,
$r>0$.

Theorem \ref{zeroone} and Corollary \ref{weaker2} are proved in
Section~3.

The problem of ``factorisatio numerorum'' can be put in the
framework of enumeration of multisets.  The following
description of factorisatio numerorum is taken from
\cite{KKWa}. An {\em (additive) arithmetical semigroup } is
a free commutative semigroup $G$ with identity element $1$,
generated by a countable set $P$ of ``prime'' elements,
and equipped with  an integer-valued ``degree''mapping
$\partial$ such that
\begin{itemize}
\item[({\em i})] $\partial(1)=0$, $\partial(p)>0$ for all $p\in P$.
\item[({\em ii})] $\partial(ab)=\partial(a)+\partial(b)$ for all $a,b\in G$.
\item[({\em iii})] The number $G^\#(j)$ of primes of degree $j$ in
$G$ is finite for all integers $j$.
\end{itemize}

Multisets can be put into the framework of arithmetical semigroups
by letting the operator $\partial$ stand for the size of the
multiset and defining the product of two multisets to be their
disjoint union. The identity element 1 is then just the empty
multiset with total size  0.

Let $f(n)$ be the total number of unordered factorizations
of elements $g\in G$  with $\partial(g)=n$. Then
\cite{KKWa} shows that $$ \sum_{n=0}^\infty f(n)x^n =
\prod_{j=1}^\infty(1-x^j)^{-G^\#(j)}, \quad \vert x\vert
<1. $$ This equation is just (\ref{ogf}), except that $c_n$
has been replaced by $f(n)$ and $a_j$  by $ G^\#(j)$. A
typical example considered in \cite{KKWa} is polynomials
over finite fields, for which $G^\#(j)=q^j$ for some prime
$q>1$.

We are able to extend the results of \cite{KKWa} and give
asymptotic results for ``factorisatio numerorum'' when
$a_j=G^\#(j)=Kj^{r-1}y^j+O(y^{\nu j})$ for $K>0,\ \nu\in(0,1)$ and
$r>0,\ y> 1$. This involves getting precise enough estimates of
$\sigma_n$ in order to derive first order asymptotics of $c_n$. We
restrict to the case $y>1$, as then a fairly simple argument using
the Poisson summation formula is effective.
\begin{theorem}\label{firstorder}
Assume that $a_j=Kj^{r-1}y^j+O\left(y^{\nu j}\right)$, where
$K>0, r>0$, $y>1$, and $\nu\in (0,1)$. Then $c_n$ has
asymptotics $$ c_n\sim \kappa_1y^nn^{-(r+2)/2(r+1)}
\exp\left(\kappa_2n^{r/(r+1)}\right)$$ for positive constants
$\kappa_1$, $\kappa_2$. Moreover, \eq\label{kappadef}
\kappa_2=\frac{r+1}{r}(K\Gamma(r+1))^{1/(r+1)},
\en
where $\Gamma$ is the gamma function.

\end{theorem}

For $r=1,$ (\ref{kappadef}) recovers the asymptotic
formula in \cite{KKWa}.

Theorem \ref{firstorder} is proved in Section 4.

\begin{remark}
It is known (\cite{Bur}, p.34) that the
generating function $g$ for $c_n$ can be written as
\eq\label{star} g(x)=\exp{S^*(x)}, \quad \vert x\vert\le 1, \en
where
$$S^*(x)=\sum_{j=1}^\infty a_j^* x^j, \quad \vert x\vert\le 1$$ is
the so-called star transformation of the generating function
$$S(x)=\sum_{j=1}^\infty a_jx^j, \ \vert x\vert \le 1$$ for $a_j,$
namely

$$a_j^*=\sum_{lk=j}\frac{a_l}{k}, \quad j\ge 1.$$

The representation (\ref{star}) says (see e.g. \cite{FG2})  that
$g$ can be viewed also as a generating function for the
parameters $a_j^*=\frac{m_j^*}{j!}, \ j\ge 1$ of  assemblies. By
Lemma 5.2 in \cite{BB}, we have that the asymptotic formula in
Theorem \ref{firstorder} for enumeration of expansive multisets
with $y>1$ is also valid for enumeration of assemblies with the
same parameters $a_j$. In this connection observe that, under
the assumption (\ref{exp}), the orders of the quantities
$\delta_n, \ B^2_n, \ \rho_l(n)$ found in Section 2  appear to be
the same as the ones in \cite{FG2} for expansive assemblies.
Summing this up, we see that the asymptotic behavior of expansive
assemblies and multisets is alike. We will show further on that
the same is true also for selections.
\end{remark}


We define the selections determined  by the parameters $a_j$ to be
those multisets for which no component type appears more than
once. For example, if $a_j=1$ for all $j$ then a selection is
an integer partition with distinct parts. Let $\tilde c_n$
denote the number of selections of size $n$ determined by the
$a_j$. Then the generating function $\tilde g$ for the $\tilde
c_n$ is

$$\tilde g(x)=\prod_{j=1}^\infty(1+x^j)^{a_j},\quad \vert x\vert\le 1$$
 and analysis similar to the one for multisets
gives that in this case $j^{-1}X_j$ is a binomial r.v. with
parameters $a_j$ and  $
\frac{\exp{(-j\sigma)}}{1+\exp{(-j\sigma)}}$, where $\sigma>0$
is arbitrary. Consequently, we have
\begin{theorem}\label{selection}
Assume that $a_j$  satisfy (\ref{exp}). Let
$\tilde\sigma_n$ be chosen in such a way that $$ \tilde
M_n:=\sum_{j=1}^n\frac{ja_je^{-j\tilde
\sigma_n}}{1+e^{-j\tilde\sigma_n}}=n.
$$
and define
 $\tilde B_n$ by \eq\label{Bsdef} (\tilde B_n)^2=
\sum_{j=1}^n\frac{j^2a_je^{-j\tilde
\sigma_n}}{(1+e^{-j\tilde\sigma_n})^2}. \en Then the number of
selections is asymptotically \eq\label{csasymp} \tilde
c_n\sim\frac{e^{n\tilde \sigma_n}}{\sqrt{2\pi (\tilde B_n)^2}} \,
\prod_{j=1}^n \left( {1+e^{-j\tilde \sigma_n}}\right)^{a_j}. \en
Moreover, if we assume that $a_j$ is as in Theorem
\ref{firstorder} then $\tilde c_n$ has the same asymptotics as
$c_n$, with a different constant $k_1$.
\end{theorem}

We sketch the proof of Theorem \ref{selection}, which is similar
to the proof of Theorems \ref{main} and Theorem \ref{firstorder},
in Section 5.

The classic example of an expansive multiset is integer
partitions. For partitions of an integer, $a_j=1$ for all
$j$, so that $r=y=1$. We can derive the Hardy-Ramanujan
formula giving asymptotics of $c_n$ for partitions from
Theorem \ref{main} by using well known properties of the
Euler generating function
$F(x):=\prod_{j=1}^\infty\left(1-x^{j}\right)^{-1},\ \
\vert x\vert < 1$. Since
$\int_0^\infty\frac{ue^{-u}}{1-e^{-u}}du=\frac{\pi^2}{6}$
(see \cite{GR}, Formula 3.411-7), we apply the
Euler-Maclaurin formula (described in detail in \cite{GKP})
to obtain from (\ref{12})
$$n=\frac{\pi^2}{6\sigma_n^2}-{1\over {2\sigma_n}} + O(1).$$
Consequently, \eq\label{sigmaasymp}
\sigma_n=\frac{\pi}{\sqrt{6n}} - \frac{1}{4n} +
O(n^{-3/2}). \en The equality (8.6.1) in Section 8.6 of
\cite{H} gives \eq\label{Fasymp} F(e^{-\sigma_n})\sim
\left(\frac{\sigma_n}{2\pi}\right)^{1/2}
\exp\left(\frac{\pi^2}{6\sigma_n}\right). \en Finally,
Theorem \ref{main} and (\ref{Fasymp}) produce
\eq\label{casymp1} c_n\sim \frac{e^{n\sigma_n}}{\sqrt{2\pi
B_n^2}} \left(\frac{\sigma_n}{2\pi}\right)^{1/2}
\exp\left(\frac{\pi^2}{6\sigma_n}\right). \en An asymptotic
analysis using (\ref{sigmaasymp}) which we do not present
here shows that the asymptotic relation (\ref{casymp1})
still holds when $\sigma_n$ is replaced by $\pi/\sqrt{6n}$.
Furthermore, $$ \sigma_n^3B_n^2\sim
\int_0^\infty\frac{u^2e^{-u}}{(1-e^{-u})^2}du
=\frac{\pi^2}{3}
$$ (see \cite{GR}, Formula 3.423-3). The above analyis results
in the Hardy-Ramanujan formula $$
c_n\sim\frac{e^{C\sqrt{n}}}{4n\sqrt{3}} \ \ \ {\rm with} \ \ \
C=\pi(2/3)^{1/2}. $$

Examples of expansive multisets, most of them with
$r=1,$ can be found in \cite{BB},\cite{Bur}. The
simplest example is the class of finite $k$-colored linear
forests, which has $a_j=k^j,$ so that $r=1$, $y=k$. We
give an example with $r=1/2$, $y=2$.

\begin{example}
Consider the linear forests in which every tree
on $j$ vertices is 2-colored with colors red and
blue in such a way that it has  exactly $[j/2]$
red vertices and $j-[j/2]$ blue vertices. Then $a_j={j\choose
{\lfloor j/2\rfloor}}\sim \sqrt{2/\pi}2^jj^{-1/2}$.
\end{example}

We may generalize the last example to get multisets with any $r\in [1/2,1]$
as follows.

\begin{example}
Consider the forests in which every component of size $j$ is
composed of a linear tree on $\lfloor j^\alpha\rfloor,$ $\alpha\in
[0,1]$ vertices and a cycle on the $j-\lfloor j^\alpha\rfloor+1$
vertices, in such a way that  one end vertex of the tree is
identified with one vertex of the cycle. Call such
components(=graphs) lollipops. Suppose we are considering
2-colored lollipops, such that the number of red vertices in the
tree is $[j^\alpha/2]$ and the number of blue vertices is
$j^\alpha-[j^\alpha/2]$. There is no restriction on the number of
blue/red vertices in the cycle. The number  of 2-colored lollipops
is $$ {{\lfloor j^\alpha\rfloor}\choose {\lfloor\lfloor
j^\alpha\rfloor/2\rfloor}}2^{j-\lfloor j^\alpha\rfloor}\sim
\sqrt{2/\pi}\, 2^{\lfloor j^\alpha\rfloor}/\sqrt{\lfloor
j^\alpha\rfloor}\times 2^{j-\lfloor j^\alpha\rfloor}
\sim\sqrt{2/\pi}\,2^jj^{-\alpha/2}. $$ This example has
$r=1-\alpha/2, \ y=2$.
\end{example}

 The next example is a natural case where the multiset
satisfies Theorem~\ref{zeroone} but not the conditions in
\cite{Bl,BB}.

\begin{example}
Consider linear forests which are $k$
colored. If a tree is $j$ vertices long and $j$ is even,
then it may be $k$-colored in all $k^j$ possible ways. If
$j$ is odd, then the first vertex is always red and the
remaining $j-1$ vertices may be colored in all $k^{j-1}$
possible ways. Then $a_j=k^j$ if $j$ is even and
$a_j=k^{j-1}$ if $j$ is odd. Therefore, (\ref{exp})
holds with $r=1,y=k, D_1=k^{-1}, D_2=1.$
\end{example}

\begin{example}
Finally, note that $r=2$, $y=1$ corresponds to plane
partitions; see \cite{Andrews}.
\end{example}

\section{Asymptotics for expansive multisets}
In this section we prove Theorem \ref{main} and Corollary
\ref{weaker1}. Recall that we assume
 here, with the exception of Corollary \ref{weaker1},
that $a_j$ obey the condition (\ref{exp}).  We first derive
an expansion for the characteristic function $\phi$
given by (\ref{phidef}) to general precision. For any
$\sigma>0,$ we define the quantities $\rho_l=\rho_l(n)$ for
$l\geq 3$ by \eq\label{rhodef} \rho_l:=
\sum_{j=1}^nj^la_j\sum_{k=1}^\infty
k^{l-1}e^{-jk\sigma}.\en
\begin{lemma}\label{phiexp}
For a fixed $n$ and any integer $s\geq 3,$  the function
$\phi$ can be expanded as $$ \phi(\alpha)= \exp\left(2\pi i
M_n\alpha - 2\pi^2 B_n^2 \alpha^2
+\sum_{l=3}^{s-1}\frac{(2\pi i)^l\rho_l}{l!}\alpha^l+
O(\alpha^s\rho_s)\right),\quad \alpha\to 0, $$

where $M_n$ and $B_n$ are given by (\ref{Mdef}) and
(\ref{Bdef}).

\end{lemma}
\begin{proof}
 The definition (\ref{phikdef}) implies that
for all $\alpha \in R$
\begin{eqnarray*}
\phi(\alpha)&=& \prod_{j=1}^n\left(
\frac{1-e^{-j\sigma}}{1-e^{-j\sigma_n}e^{2\pi i \alpha j}}
\right)^{a_j}\\ &=& \exp\left(\sum_{j=1}^na_j\Bigg(
\log\left(1-e^{-j\sigma}\right)-
\log\left(1-e^{-j\sigma}e^{2\pi i \alpha j}\right)
\Bigg)\right).
\end{eqnarray*}
The logarithms may be expanded in Taylor series as
$\sigma>0$ and $\alpha\in R$ are fixed, giving
\begin{eqnarray*}
\phi(\alpha) &=& \exp\left(\sum_{j=1}^na_j
\left(-\sum_{k=1}^\infty\frac{e^{-\sigma jk}}{k}+
\sum_{k=1}^\infty\frac{e^{-\sigma jk}e^{2\pi ijk\alpha}}{k}
\right)\right).
\end{eqnarray*}
We make use of the Taylor expansion with $s\ge 3$
$$e^{2\pi ijk \alpha}=1+2\pi ijk\alpha - 2\pi^2j^2k^2\alpha^2 +
\sum_{l=3}^{s-1} \frac{(2\pi ijk\alpha)^l}{l!}+ O(j^sk^s\alpha^s),
\quad \alpha\to 0, $$ which holds uniformly for all $j,k\geq 1$,
to get
\begin{eqnarray} \label{chf}
\phi(\alpha) &=& \exp\left(\sum_{j=1}^na_j \left(
\sum_{k=1}^\infty\frac{e^{-jk\sigma}}{k} \left[2\pi
ijk\alpha-2\pi^2j^2k^2\alpha^2 +\sum_{l=3}^{s-1}
\frac{(2\pi ijk\alpha)^l}{l!} +O\left(\alpha^sj^sk^s\right)
\right] \right)\right)\nonumber\\ &=& \exp\left( 2\pi
i\sum_{j=1}^n
\frac{ja_je^{-j\sigma}}{1-e^{-j\sigma_n}}\alpha
-2\pi^2\sum_{j=1}^n\frac{j^2a_je^{-j\sigma}}{(1-e^{-j\sigma_n})^2}\alpha^2
+\sum_{l=3}^{s-1}\frac{(2\pi i)^l\rho_l}{l!}\alpha^l
+O(\rho_s\alpha^s)\right)\nonumber\\ &=& \exp\left(2\pi i
M_n\alpha - 2\pi^2 B_n^2 \alpha^2
+\sum_{l=3}^{s-1}\frac{(2\pi i)^l\rho_l}{l!}\alpha^l+
O(\alpha^s\rho_s)\right), \ \alpha\to 0.
\end{eqnarray}
\end{proof}

In what follows we set $\sigma=\sigma_n$ determined by
(\ref{12}) and define $\delta_n$ by
$\delta_n:=\sigma_n-\log y, \ y\ge 1$. In proving Theorem
\ref{main} we will apply Lemma \ref{phiexp} and so we need
estimates of $\delta_n$, $B_n$ and $\rho_3$. (We use the
ability to expand $\phi(\alpha)$ to higher order precision
than $s=3$ for the proof of Theorem \ref{zeroone} in
Section 3.)

\begin{lemma}\label{asympsigma}
$\delta_n\asymp n^{-1/(r+1)}$, $B_n^2\asymp n^{(r+2)/(r+1)}$,
and $\rho_l(n)\asymp n^{(r+l)/(r+1)}$, $r>0$, for all $l\geq
3$. Moreover, there exists $N>0$  such that
$\delta_{n+1}<\delta_n$, $B_n^2\leq B_{n+1}^2$, and
$\rho_l(n)\leq\rho_l(n+1)$ whenever $n\geq N.$
\end{lemma}

\begin{proof}
We first prove some preliminary facts about $\delta_n$. Let
$D_1,D_2>0$ be constants such that $D_1j^{r-1}y^j\leq a_j\leq
D_2j^{r-1}y^j, \quad j\ge 1,\quad y\ge 1,quad r>0$. Since
$\sigma_n>0, \ n\ge 1$, by (\ref{12}),
 we deduce from
(\ref{Mdef}),(\ref{12}) that $n\geq\sum_{j=1}^n j a_j
e^{-j\sigma_n}\geq D_1\sum_{j=1}^nj^re^{-j\delta_n}, $
implying that $\delta_n>0$ for $n$ large enough. Suppose
that there exists a constant $\epsilon>0$ and a subsequence
$n_k\to\infty$ such that $\delta_{n_k}\geq \epsilon$. Then,
again by (\ref{Mdef}),(\ref{12}), $$ n_k\leq \frac{
D_2}{1-y^{-1}e^{-\epsilon}}
\sum_{j=1}^{n_k}j^re^{-j\epsilon}=O(1),\quad n_k\to \infty.
$$ We therefore must have $\delta_n\to 0$.
 Next, we derive from the inequality
$$n\ge D_1\sum_{j=1}^{n}j^re^{-j\delta_{n}}\geq
D_1e^{-n\delta_{n}}\sum_{j=1}^{n} j^r, \ r>0 $$ that
$n\delta_n\to\infty$ as $n\to\infty$.

With the help of these facts we further get $$ n\leq
D_2\sum_{j=1}^{n}\frac{j^re^{-j
\delta_{n}}}{1-y^{-j}e^{-j\delta_{n}}} \leq
\delta_{n}^{-r-1}D_2\sum_{j=1}^{\infty}\frac{(j\delta_{n})^re^{-j\delta_{n}}}
{1-e^{-j\delta_{n}}}\delta_{n} \sim
\delta_{n}^{-r-1}D_2\int_0^\infty\frac{x^re^{-x}}{1-e^{-x}}dx.
$$

Since the last  integral is  bounded,  we conclude that
$\delta_n\le D_4 n^{-1/(r+1)}, n\ge 1,$  where $D_4>0$ is a
constant. On the other hand,

$$ n\geq D_1\sum_{j=1}^{n} j^re^{-j\delta_{n}} \sim
\delta_{n}^{-r-1}D_1\int_0^\infty x^re^{-x}dx. $$ This gives
 $\delta_n\geq
D_3 n^{-1/(r+1)}, n\ge 1,$  where $D_3>0$ is a constant.
 We have shown
that $\delta_n\asymp n^{-1/(r+1)}$.

For any $l\geq s\geq 0$, arguments similar to those above show
that \eq\label{11}
\sum_{j=1}^n\frac{j^la_je^{-j\sigma_n}}{(1-e^{-j\sigma_n})^s}
\asymp
\sum_{j=1}^n\frac{j^{r+l-1}e^{-j\delta_n}}{(1-y^{-j}e^{-j\delta_n})^s}
\asymp \delta_n^{-r-l}\asymp n^{(r+l)/(r+1)}.\en The last
asymptotic applied to  (\ref{Bdef}) results in the stated
asymptotics for $B_n^2$.

We now show  that $\rho_l, \ l\ge 3$ has the same
asymptotics as in (\ref{11}). We have
\begin{eqnarray*} \label{25}
\rho_l\asymp \sum_{j=1}^nj^{r+l-1}\sum_{k=1}^\infty
k^{l-1}e^{-\delta_n jk}y^{j-jk}\\ \sim
\delta_n^{-r-l}\sum_{k=1}^\infty\Big(
k^{-r-1}\int_{k\delta_n}^{kn\delta_n} z^{r+l-1}
y^{\frac{z}{k\delta_n}{(1-k)}}e^{-z}dz\Big), \quad r\ge 0.
\end{eqnarray*}
The integral in the last expression is $\le \Gamma(r+l)$
for all $k\ge 1,$ with equality for $k=1,$
 since
$n\delta_n\to \infty, \ \delta_n\to 0, \ n\to \infty, $
while $y^{\frac{z}{k\delta_n}{(1-k)}}\le 1, \ k\ge 1, \
y\ge 1, z\ge 0.$ Thus, the last series in (\ref{11})
converges which implies
 \eq \label{13} \rho_l\asymp\delta_n^{-r-l}\asymp
n^{(r+l)/(r+1)}.\en

Lastly we will prove that $\sigma_n$, and so $\delta_n$, is
eventually monotone decreasing in $n$. Suppose that
$\sigma_{n+1}\geq \sigma_n$ for some $n$. Then, it
would follow that $$
\frac{e^{-j\sigma_{n+1}}}{1-e^{-j\sigma_{n+1}}}\leq
\frac{e^{-j\sigma_n}}{1-e^{-j\sigma_n}} $$ for all $j\ge
1$, and consequently
\begin{eqnarray*}
n+1&=&\sum_{j=1}^{n+1}\frac{ja_je^{-j\sigma_{n+1}}}{1-e^{-j\sigma_{n+1}}}
\leq\sum_{j=1}^{n+1}\frac{ja_je^{-j\sigma_n}}{1-e^{-j\sigma_n}}
=
n+\frac{(n+1)a_{n+1}y^{-n-1}e^{-(n+1)\delta_n}}{1-y^{-n-1}e^{-(n+1)\delta_n}}.
\end{eqnarray*}

In view of the established asymptotics for $\delta_n$, the  last
term  in the preceding inequality tends to $0$ for all $y\ge 1.$
We therefore must have $\sigma_{n+1}<\sigma_n$ for
sufficiently large $n$ which implies  the same inequality
for $\delta_n.$

The derivative $\frac{d}{dx}\left[e^{-x}/(1-e^{-x})^2\right]$ is
negative for $x\geq 0$, so that if $n$ is sufficiently large,

 $$
B_n^2=\sum_{j=1}^n\frac{j^2a_je^{-j\sigma_n}}{(1-e^{-j\sigma_n})^2}
<\sum_{j=1}^{n+1}\frac{j^2a_je^{-j\sigma_{n+1}}}{(1-e^{-j\sigma_{n+1}})^2}
=B_{n+1}^2 $$ and similarly,  $\rho_l(n)\leq\rho_l(n+1)$
when $n$ is large.
\end{proof}

Lemma \ref{local} below proves a local limit theorem for the
probability in (\ref{zz}). Define the sequence $\alpha_0(n)$
by \eq\label{alpha0def} \alpha_0(n):=\delta_n^{(r+2)/2}\log
n\asymp n^{-(r+2)/2(r+1)}\log n. \en We will express the integral
in (\ref{zz}) as
 \eq\label{Tdef} T= T(n):= \int_0^1 \phi(\alpha)\,
e^{-2\pi i\alpha n} d\alpha = \int_{-1/2}^{1/2}
\phi(\alpha)\, e^{-2\pi i\alpha n} d\alpha = T_1(n)+T_2(n),
\en where the middle equality  follows from the
periodicity of $\phi(\alpha)$, as defined by
(\ref{phidef}), (\ref{phikdef}), and where
\eq\label{T1def}
T_1=T_1(n)=\int_{-\alpha_0(n)}^{\alpha_0(n)} \phi(\alpha)\,
e^{-2\pi i\alpha n}d\alpha, \en
\eq\label{T2def}T_2=T_2(n)=\int_{-1/2}^{-\alpha_0(n)}
\phi(\alpha)\, e^{-2\pi i\alpha n} d\alpha +
\int_{\alpha_0(n)}^{1/2} \phi(\alpha)\, e^{-2\pi i\alpha n}
d\alpha. \en

Lemma \ref{local} and the representation
(\ref{rewrite1}) prove Theorem \ref{main}.
\begin{lemma}\label{local}
\eq\label{T1est} T_1\sim \left(2\pi
B_n^2\right)^{-1/2}
\en
and for sufficiently large $n$,\eq\label{T2est}
T_2\leq\exp\left(-C\log^2 n\right),
\en
for a constant $C>0$, from which it follows that
$$
T=\int_0^1 \phi(\alpha)e^{-2\pi i\alpha n}\, d\alpha
\sim\left(2\pi B_n^2\right)^{-1/2}.
$$
\end{lemma}
\begin{proof}
The proof of this lemma is similar to the proof of
Lemma 6 in \cite{FG}. Using the expansion of Lemma
\ref{phiexp} with $s=3$ in the definition (\ref{T1def}) of
$T_1$ and observing that, by virtue of Lemma
\ref{asympsigma}, $\lim_{n\to\infty}\alpha^3\rho_3=0$ and
$\alpha^2 B_n^2\to \infty$ for all
$\alpha\in[-\alpha_0(n),\alpha_0(n)],$ gives
$$ T_1\sim\int_{-\alpha_0(n)}^{\alpha_0(n)} \exp\left(-
2\pi^2\alpha^2 B_n^2\right)d\alpha \sim \left(2\pi
B_n^2\right)^{-1/2}.
$$

The bound for $T_2$ starts with the identity for
all $\alpha\in R$,
$$ |\phi(\alpha)|=
\prod_{j=1}^n\left|\frac{1-e^{-j\sigma_n}}{1-e^{-
j\sigma_n}e^{2\pi i \alpha j}} \right|^{a_j}\\
=
\exp\left(-\sum_{j=1}^n {a_j\over
2}\log\left(1+\frac{4e^{-j\sigma_n}\sin^2(\pi\alpha j)}
{\left(1-e^{-j\sigma_n}\right)^2}\right)\right). $$ All of the
logarithms are positive, and $\log(1+x)\geq x/(1+c)$ whenever
$0\le x\le c$ for a constant $c>0$, so we have
\begin{eqnarray*}
|\phi(\alpha)|&\leq& \exp\left(-\sum_{(4\sigma_n)^{-1}\leq j\leq
n} {a_j\over 2}\log\left(1+\frac{4e^{-j\sigma_n}\sin^2(\pi\alpha
j)} {\left(1-e^{-j\sigma_n}\right)^2}\right)\right)\\ &\leq&
\exp\left(-\sum_{(4\sigma_n)^{-1}\leq j\leq n}
C_1a_je^{-j\sigma_n}\sin^2(\pi\alpha j) \right)\\ &\leq&
\exp\left(-\sum_{(4\sigma_n)^{-1}\leq j\leq n}
C_2j^{r-1}e^{-j\delta_n}\sin^2(\pi\alpha j) \right), \ \alpha\in
R,
\\
\end{eqnarray*}
for some constants $C_1,C_2>0$. In view of the inequality $\delta_n\leq \delta_n +\log y=\sigma_n, \ y\ge 1$, we
also have \eq\label{phibasic} |\phi(\alpha)|\leq \exp\left(-\sum_{(4\delta_n)^{-1}\leq j\leq n}
C_2j^{r-1}e^{-j\delta_n}\sin^2(\pi\alpha j) \right):=\exp (-V_n(\alpha)), \ \alpha\in R.\en

Since our $\delta_n$ is of the same order as $\sigma_n$ in \cite{FG},
the argument of Lemma 7 in
\cite{FG}, gives the desired estimate of $ V_n(\alpha)$:
 \eq\label{Vest}
V_n(\alpha)\asymp \delta_n^{-r}\asymp n^{r/(r+1)}\gg \log^2 n, \ \alpha\in [\alpha_0, 1/2]
\en
\end{proof}


\begin{proof}[of Corollary \ref{weaker1}]
 For $a\in\calF(2r/3+\epsilon,r,y)$, $r>0$, $0<\epsilon\le r/3$,
$y\geq 1$, the arguments in Lemma \ref{asympsigma} show
that \eq\label{asymp1} n^{-1/(r_1+1)}\preceq \delta_n
\preceq n^{-1/(r_2+1)},  {\rm \ \ \ \ \ \ }
\delta_n^{-(r_1+2)}\preceq B_n^2\preceq\delta_n^{-(r_2+2)},
\en \eq\label{asymp2} \delta_n^{-(r_1+l)}\preceq
\rho_l(n)\preceq\delta_n^{-(r_2+l)} {\rm \ \ for \ }l\geq
3. \en We write $\alpha_0=(B_n)^{-1}\log n$ to obtain, as
$n\to \infty,$ $$ \alpha_0^3\rho_3(n)\leq\gamma_3(\log^3
n)\delta_n^{3(r_1+2)/2} \delta_n^{-(r_2+3)}=\gamma_3(\log^3
n)\delta_n^{(3r_1-2r_2)/2}\to 0,$$ since $3r_1-2r_2>0$. We
have shown that the asymptotic (\ref{T1est}) is still
valid. The upper bounds on $T_2$ are like those in the
proof of  Theorem \ref{main}, with the replacement of
(\ref{Vest}) by
$$\delta_n^{-r_1}\preceq V_n(\alpha)\preceq
\delta_n^{-r_2} {\rm \ \ for \ } \alpha\in [\alpha_0,1/2].
$$
\end{proof}

\section{Logical limit laws for expansive multisets}
Lemma \ref{ratio} below and the asymptotic
$\delta_{n+1}\asymp n^{-1/(r+1)}$ from Lemma \ref{asympsigma} show
that the $c_n$ satisfy the hypotheses of Theorem \ref{basic} and
thereby prove Theorem \ref{zeroone}.
\begin{lemma}\label{ratio}
If $a_j\asymp j^{r-1}y^j$, with $r>0$ and $y\geq 1$, then $$
\frac{c_n}{c_{n+1}}=y^{-1}e^{-\delta_n+o\left(\delta_n\right)}.
  $$
\end{lemma}
\begin{proof}
We use (\ref{rewrite1}) to get \eq\label{ratio1}
\frac{c_n}{c_{n+1}}=
e^{n\sigma_n-(n+1)\sigma_{n+1}}\prod_{j=1}^n \left(
\frac{1-e^{-j\sigma_{n+1}}} {1-e^{-j\sigma_n}}
\right)^{a_j}
\left(1-e^{-(n+1)\sigma_{n+1}}\right)^{a_{n+1}}
\frac{T(n)}{T(n+1)}, \en where $T(n)$ is defined by
(\ref{Tdef}). Since $$ e^{-(n+1)\sigma_{n+1}}a_{n+1}\asymp
(n+1)^{r-1}e^{-(n+1)\delta_{n+1}}, $$ by definition of
$\delta_n$, and $$ (n+1)\delta_{n+1}\asymp (n+1)^{r/(r+1)},
$$ by Lemma \ref{asympsigma}, it follows that
\eq\label{termsmall1}
\left(1-e^{-(n+1)\sigma_{n+1}}\right)^{a_{n+1}}=e^{o(\delta_n)}.
\en The second factor in the RHS of (\ref{ratio1}) may
be rewritten as \eq\label{ratio2} \prod_{j=1}^n \left(
\frac{1-e^{-j\sigma_{n+1}}} {1-e^{-j\sigma_n}}
\right)^{a_j} = \exp\left( \sum_{j=1}^n a_j \log
\left(1-\frac{e^{-j\sigma_{n+1}}-e^{-j\sigma_n}}
{1-e^{-j\sigma_n}} \right)\right).\en We assume that $n>N$
for the $N$ in Lemma \ref{asympsigma}, so that in
particular $\sigma_{n+1}<\sigma_n$. Since $\log(1-x)\leq
-x$ when $x\in [0,1]$, we have
\begin{eqnarray*}
\prod_{j=1}^n
\left(
\frac{1-e^{-j\sigma_{n+1}}}
{1-e^{-j\sigma_n}}
\right)^{a_j}
&\leq&
\exp\left(-\sum_{j=1}^na_j\frac{e^{-j\sigma_{n+1}}-e^{-j\sigma_n}}
{1-e^{-j\sigma_n}}\right)\\
&\leq&
\exp\left(-\sum_{j=1}^na_j\frac{e^{-j\sigma_n}(j\sigma_n-j\sigma_{n+1})}
{1-e^{-j\sigma_n}}\right)\\
&=&
\exp\left(-(\sigma_n -\sigma_{n+1})\sum_{j=1}^n\frac{ja_je^{-j\sigma_n}}
{1-e^{-j\sigma_n}}\right)\\
&=&
 \exp\left(-(\sigma_n -\sigma_{n+1})M_n
\right)\\
&=&
e^{-(\sigma_n-\sigma_{n+1})n},
\end{eqnarray*}
where the second inequality results from the fact that
$e^z-1\ge z, \ z\ge 0$ . Since $\log(1-x)\geq -x/(1-x)$ for
$x\in [0,1]$,
we lower bound (\ref{ratio2}) by
\begin{eqnarray*}
\prod_{j=1}^n
\left(
\frac{1-e^{-j\sigma_{n+1}}}
{1-e^{-j\sigma_n}}
\right)^{a_j}
&\geq&
\exp\left(-\sum_{j=1}^na_j\frac{e^{-j\sigma_{n+1}}-e^{-j\sigma_n}}
{1-e^{-j\sigma_{n+1}}}\right)\\
&\geq&
\exp\left(-\sum_{j=1}^na_j\frac{e^{-j\sigma_{n+1}}(j\sigma_n-j\sigma_{n+1})}
{1-e^{-j\sigma_{n+1}}}\right)\\
&=& \exp\left(-(\sigma_n-\sigma_{n+1})\sum_{j=1}^{n+1}
\frac{ja_je^{-j\sigma_{n+1}}}
{1-e^{-j\sigma_{n+1}}}\right)\\
&\ge& \exp\left(-(\sigma_n -\sigma_{n+1})M_{n+1}
\right)\\
&=&
e^{-(\sigma_n-\sigma_{n+1})(n+1)}.
\end{eqnarray*}
Thus, the product of the first two factors of
(\ref{ratio1}) is bounded above and below by
$$
e^{-\sigma_n}\leq e^{n\sigma_n-(n+1)\sigma_{n+1}}\prod_{j=1}^n \left( \frac{1-e^{-j\sigma_{n+1}}}
{1-e^{-j\sigma_n}} \right)^{a_j} \leq e^{-\sigma_{n+1}}=e^{-\sigma_n +(\sigma_n-\sigma_{n+1})}.
$$
We bound $\sigma_n-\sigma_{n+1}$ by observing that
\begin{eqnarray*}
1&=&(n+1)-n\\
&=&\sum_{j=1}^{n+1}
\frac{ja_je^{-j\sigma_{n+1}}}
{1-e^{-j\sigma_{n+1}}}-\sum_{j=1}^n
\frac{ja_je^{-j\sigma_n}}
{1-e^{-j\sigma_n}}\\
&\geq&
\sum_{j=1}^n
\frac{ja_j(e^{-j\sigma_{n+1}}-e^{-j\sigma_n})}
{1-e^{-j\sigma_n}}\\
&\geq& (\sigma_n-\sigma_{n+1}) \sum_{j=1}^n
\frac{j^2a_je^{-j\sigma_n}} {1-e^{-j\sigma_n}},
\end{eqnarray*}
for $n$ sufficiently large.

Thus,  recalling that $\sigma_n-\sigma_{n+1}>0$, it follows from
(\ref{11}), applied with   $l=2$ and $s=1$, and Lemma
\ref{asympsigma} that \eq\label{deldiff}
\sigma_n-\sigma_{n+1}\le O\left(\delta_n^{r+2}\right).
\en We have shown that
\begin{eqnarray}\label{termsmall2}
e^{n\sigma_n-(n+1)\sigma_{n+1}}
\prod_{j=1}^n
\left(
\frac{1-e^{-j\sigma_{n+1}}}
{1-e^{-j\sigma_n}}
\right)^{a_j}
&=&\exp\left(-\sigma_n+o\left(\delta_n\right)\right)\nonumber\\
&=&y^{-1}\exp\left(-\delta_n+o\left(\delta_n\right)\right).
\end{eqnarray}

Because of (\ref{ratio1}), (\ref{termsmall1}) and
(\ref{termsmall2}), the proof will be completed if we show
that $\frac{T(n)}{T(n+1)}= e^{o\left(\delta_n\right)}.$ The
definitions (\ref{Tdef}), (\ref{T1def}), (\ref{T2def})
along with Lemma \ref{asympsigma} and Lemma \ref{local}
imply that $$ \frac{T_2(n)}{T_1(n)}=o(\delta_n),$$ which
gives
\begin{eqnarray}\label{Tratio}
\frac{T(n)}{T(n+1)} &=&1+\frac{T_1(n)-T_1(n+1)}{T_1(n+1)} +
o\left(\delta_n\right).
\end{eqnarray}
The definition (\ref{T1def}) together with (\ref{alpha0def}) and
(\ref{deldiff}) produce $$ \vert T_1(n)-T_1(n+1)\vert \le
\int_{-\alpha_0(n)}^{\alpha_0(n)} \big\vert
\phi_n(\alpha)e^{-2\pi i \alpha n}-\phi_{n+1}(\alpha)e^{-2\pi i
\alpha (n+1)}\big \vert \,d\alpha + o(\delta_n).$$

Next, Lemma
\ref{phiexp}, (\ref{12}) and the monotonicity of $\rho_l(n)$ imply that, for
a sufficiently large
fixed $n$ and $\alpha\to 0,$
\begin{eqnarray}
\label{bar} &&
\phi_n(\alpha)e^{-2\pi i \alpha n}-\phi_{n+1}(\alpha)e^{-2\pi i \alpha (n+1)}\nonumber\\
&=&
\exp\left(-2\pi^2\alpha^2 B_n^2
+Q_s(\alpha,n)+ O(\alpha^s\rho_s(n))\right)\\
&& - \exp\left(-2\pi^2\alpha^2 B_{n+1}^2 +Q_s(\alpha,n+1)+
O(\alpha^s\rho_s(n+1))\right)\nonumber\\&=& \phi_n(\alpha)e^{-2\pi
i \alpha n} \nonumber\\&& \times\Bigg\{1-\exp\left(-2\pi^2\alpha^2
(B_{n+1}^2-B_n^2) + Q_s(\alpha,n+1)-Q_s(\alpha,n)+
O(\alpha^s\rho_s(n+1))\right)\Bigg\}\nonumber,
\end{eqnarray}
where we denoted
$$Q_s(\alpha,n)=\sum_{l=3}^{s-1}(2\pi i)^l\frac{\rho_l(n)}{l!}\alpha^l.$$

We now apply (\ref{deldiff}),  Lemma \ref{asympsigma} and
(\ref{11}) with $l=3, s=2$ to get
\begin{eqnarray}\label{Bdiff}
0\leq B_{n+1}^2-B_n^2&=& \sum_{j=1}^{n+1}
\frac{j^2a_je^{-j\sigma_{n+1}}}
{(1-e^{-j\sigma_{n+1}})^2}-\sum_{j=1}^n
\frac{j^2a_je^{-j\sigma_n}} {(1-e^{-j\sigma_n})^2}\nonumber\\
&\leq& \sum_{j=1}^n \frac{j^2a_je^{-j\sigma_{n+1}}}
{(1-e^{-j\sigma_{n+1}})^2}(j\sigma_n-j\sigma_{n+1})+
\frac{(n+1)^2a_{n+1}e^{-(n+1)\sigma_{n+1}}}
{(1-e^{-(n+1)\sigma_{n+1}})^2} \nonumber\\ &\le&
O\left(\delta_n^{-1}\right).
\end{eqnarray}

In a similar way we also have from (\ref{rhodef})
\begin{eqnarray}\label{rhodiff}
0\leq
\rho_l(n+1)-\rho_l(n)&\le&(\sigma_n-\sigma_{n+1}) \sum_{j=1}^{n}
j^{l+1}a_j\sum_{k=1}^\infty k^l e^{-jk\sigma_{n+1}}\nonumber\\
&&+\,(n+1)^{2l+r-1}e^{-\delta_n(n+1)}\nonumber\\ &\leq&O\left(
\delta_n\right)\sum_{j=1}^{n} j^{l+r}\sum_{k=1}^\infty k^l
e^{-jk\delta_{n+1}}\nonumber\\ &=&
O\left(\delta_n^{-l}\right)\sum_{j=1}^n
j^{r-1}\int_{j\delta_{n+1}}^\infty u^le^{-u}du \nonumber\\&\leq&
O\left(\delta_n^{-l}\right)\left(\sum_{j\in D_1}j^{r-1}+
\sum_{j\in D_2} j^{r-1}\int_{n^\epsilon}^\infty u^le^{-u}du
\right)\nonumber\\&\le&O\left(\delta_n^{-(l+r)}n^{\epsilon
r}\right), \quad l\ge 3, \quad \epsilon>0.
\end{eqnarray}
Here $D_1=[1, \delta_{n+1}^{-1}n^\epsilon]$ and $D_2=
[1,n]\setminus D_1$. For $\alpha\in
[-\alpha_0(n),\alpha_0(n)]$, (\ref{alpha0def}),
(\ref{Bdiff}) and Lemma \ref{asympsigma} imply
\eq\label{Bdiff2} \left|(B_{n+1}^2-B_n^2)\alpha^2\right|\le
O(\delta_n^{r+1}\log^2 n)\to  0.\en Similarly, for
$\alpha\in [-\alpha_0(n),\alpha_0(n)]$, $l\geq 3$ and all
$r>0$, \eq \label{56}
\left|(\rho_l(n+1)-\rho_l(n))\alpha^l\right|\le
O(\delta_n^{r(\frac{l}{2}-1)}n^{\epsilon r}\log^l n) \to 0,
\en for sufficiently small $\epsilon
>0$. Because of (\ref{alpha0def}) and Lemma \ref{asympsigma}, it
follows that for $\alpha\in [-\alpha_0(n),\alpha_0(n)],$
\eq\label{remainder}
\vert\alpha^l\vert\rho_l(n)=O\left(\delta_n^{r(l-2)/2}\log^l
n\right)\to 0,\quad l\ge 3. \en The above discussion
reveals the following remarkable feature of the choice
 $\alpha\in [-\alpha_0(n),\alpha_0(n)]$ in the
 expansions (\ref{bar}) and in Lemma \ref{phiexp}:
under this choice the terms with $s> 3$ can be ignored, as $n\to
\infty.$ Therefore, based on the preceding bounds, we get
\begin{eqnarray*}
&&\vert\phi_n(\alpha)e^{-2\pi i \alpha
n}-\phi_{n+1}(\alpha)e^{-2\pi i \alpha (n+1)}\vert\\
&\le&O\left(n^{-1}\log^2 n\right)\exp\left(-2\pi^2\alpha^2
B_n^2 +O\left(\delta_n^{r/2}\log ^3 n\right)\right),
\end{eqnarray*} uniformly for $\alpha\in[-\alpha_0(n),\alpha_0(n)]$.
Now it follows that
$$
\frac{\vert T_1(n)-T_1(n+1)\vert}{T_1(n+1)}\le
O(n^{-1}\log^2n),
$$
and therefore (\ref{Tratio}) gives
$\frac{T(n)}{T(n+1)}=1+O(n^{-1}\log^2
n)+o\left(\delta_n\right)=
1+o\left(\delta_n\right)=\exp\left(o\left(\delta_n\right)\right)$,
proving the lemma.
\end{proof}


\begin{proof}[of Corollary~\ref{weaker2}]
We will make use of  (\ref{asymp1}) and (\ref{asymp2}). The
argument in the proof of Lemma \ref{asympsigma} shows that
in the case considered we still have that  $\delta_n\to 0$
as $n\to\infty$ and that $\delta_n$ decreases monotonically
for large enough $n$. Let, as before, $r_1=2r/3+\epsilon, \
0<\epsilon\le r/3  $ and $r_2=r.$ Then, observing that the
bound in (\ref{deldiff}) is valid with $r$ replaced by
$r_1$, we get that the bound in (\ref{Bdiff}) becomes
$O\left(\delta_n^{r_1-r-1}\right)$. Consequently, setting
$\alpha_0$ as before gives (\ref{Bdiff2}) with $r$ replaced
by $r_1$. The left hand sides of (\ref{rhodiff}) and
(\ref{remainder}) may be bounded similarly.
\end{proof}

\section{Explicit asymptotic formulae for enumeration of expansive multisets}

In this section we will prove Theorem \ref{firstorder}, which
gives first order asymptotics for $c_n$ when $y>1, \ K,r>0$
and \eq\label{acond} a_j=Kj^{r-1}y^j+O(y^{\nu j}), \ \ \ {\rm
where} \ \  \nu\in (0,1). \en

To approximate $\sigma_n=\log y+\delta_n$ in the case considered,
it is necessary to analyze the equation

\begin{eqnarray}\label{firstest}
n&=&\sum_{j=1}^n\frac{ja_je^{-j\sigma_n}}{1-e^{-j\sigma_n}}\nonumber\\
&=&
\sum_{j=1}^n\frac{ja_jy^{-j}e^{-j\delta_n}}{1-y^{-j}e^{-j\delta_n}}\nonumber\\
&=&\sum_{j=1}^nja_jy^{-j}e^{-j\delta_n}+
O\left(\sum_{j=1}^nja_jy^{-2j}e^{-2j\delta_n}\right)\nonumber\\
&=& \sum_{j=1}^nKj^re^{-j\delta_n} + O\left(
\sum_{j=1}^nj^ry^{-(1-\nu)j}e^{-j\delta_n} \right)
+O(1)\nonumber\\
&=&
\sum_{j=1}^nKj^re^{-j\delta_n} + O(1).\nonumber
\end{eqnarray}
The Poisson summation formula as used in the proof of Lemma
4 of \cite{EG} shows that for $l>-1$,
 \eq\label{Granexp}
\sum_{j=1}^nj^le^{-j\delta_n}= \Gamma(l+1)\delta_n^{-l-1}+
C_l+ O\left(\delta_n\right),  \en where in the case $l>0$
the constant $C_l$ can be found explicitly:
$$
C_l=2\,\Gamma(l+1)(2\pi)^{-l-1}\zeta(l+1)\cos\frac{\pi(l+1)}{2}
$$
(here $\zeta(\cdot)$ is the Riemann zeta function).
The preceding estimates imply that
\eq\label{Starkexp}
n=K\delta_n^{-r-1}\Gamma(r+1)+O(1).
\en
from which it follows that
\begin{eqnarray*}
\delta_n&=&\left(
\frac{n}{K\Gamma(r+1)}+O(1)\right)^{-1/(r+1)}\\ &=&
\left(\frac{n}{K\Gamma(r+1)}\right)^{-1/(r+1)}+o(n^{-1})
\end{eqnarray*}
and that
$$
e^{n\sigma_n}\sim y^n\exp\left(n^{r/(r+1)}(K\Gamma(r+1))^{1/(r+1)}\right).
$$

The asymptotic for $B_n^2$ follows from
\begin{eqnarray*}
B_n^2&=&\sum_{j=1}^n\frac{j^2a_je^{-j\sigma_n}}{(1-e^{-j\sigma_n})^2}\\
&=&\sum_{j=1}^n\frac{Kj^{r+1}e^{-j\delta_n}}{(1-y^{-j}e^{-j\delta_n})^2}
+\sum_{j=1}^n\frac{j^{2}O(y^{-(1-\nu)j})e^{-j\delta_n}}{(1-y^{-j}e^{-j\delta_n})^2}\\
&=&\sum_{j=1}^nKj^{r+1}e^{-j\delta_n}
+O(1)\\
&\sim&
K\Gamma(r+2)\,\delta_n^{-r-2}\\
&\sim&
K^{-1/(r+1)}{\Gamma(r+1)}^{-(r+2)/(r+1)}\Gamma(r+2)n^{(r+2)/(r+1)}
\end{eqnarray*}

The second factor in (\ref{mainasymp}) may be expanded as
\begin{eqnarray}\label{expterm}
\prod_{j=1}^n \left(
{1-e^{-j\sigma_n}}\right)^{-a_j}
&=&
\exp\left(\sum_{j=1}^n
-a_j\log(1-y^{-j}e^{-j\delta_n})\right)\nonumber\\
&=&
\exp\left(\sum_{j=1}^n
\sum_{k=1}^\infty \frac{a_jy^{-jk}e^{-jk\delta_n}}{k}\right)\nonumber\\
&=&
\exp\left(\sum_{j=1}^n
a_jy^{-j}e^{-j\delta_n}
+\sum_{j=1}^n
\sum_{k=2}^\infty \frac{a_jy^{-jk}e^{-jk\delta_n}}{k}\right).
\end{eqnarray}
We use (\ref{Granexp}) and (\ref{Starkexp}) to show that the first
term in the exponential in (\ref{expterm}) equals
\begin{eqnarray*}
\sum_{j=1}^n
a_jy^{-j}e^{-j\delta_n}&=&
K\sum_{j=1}^n j^{r-1}e^{-j\delta_n} + \sum_{j=1}^n(a_j-Kj^{r-1})y^{-j}e^{-j\delta_n}\\
&=&
K\Gamma(r)\delta_n^{-r}+
KC_{r-1}+\sum_{j=1}^\infty y^{-j}(a_j-Kj^{r-1})+o(1)\\
&=& D_r n^{r/(r+1)}+KC_{r-1}+\sum_{j=1}^\infty y^{-j}
(a_j-Kj^{r-1}) + o(1),
\end{eqnarray*}
where
$$
D_r=K^{1/(r+1)}\Gamma(r) (\Gamma(r+1))^{-r/(r+1)}=
\frac{1}{r}{(K\Gamma(r+1)}^{1/(r+1)}.
$$
The second term in the exponential in (\ref{expterm}) equals
$$
\sum_{j=1}^n
\sum_{k=2}^\infty \frac{a_jy^{-jk}e^{-jk\delta_n}}{k}=
\sum_{j=1}^\infty
\sum_{k=2}^\infty \frac{a_jy^{-jk}}{k}+o(1),
$$
where the double sum on the right converges absolutely because
of (\ref{acond}).


\section{Asymptotics for expansive selections}
Let $\cs_n$ be the number of selections of total size $n$
corresponding to a given sequence $a_j$. The generating
function for $\cs_n$ is given by
\begin{eqnarray*}
\gs(x)&=&\sum_{n=0}^\infty \cs_n x^n\\ & =&
1+\sum_{n=1}^\infty\sum_{\vec{\eta} \in\Omega_n}
\prod_{j=1}^n{a_j\choose \eta_j}x^{j\eta_j}\\ &=&
\prod_{j=1}^\infty\left(1+x^j\right)^{a_j},\vert
x\vert <1.
\end{eqnarray*}
By adapting the derivation of (\ref{ident1}) for multisets to
the truncated generating function
$\gs_n(x)=\prod_{j=1}^n\left(1+x^j\right)^{a_j}$ we obtain
for all $\sigma\in R$, \eq\label{ident2}
\cs_n=e^{n\sigma} \int_0^1 \prod_{j=1}^n \left(
{1+e^{-j\sigma}e^{2\pi i\alpha j}}\right)^{a_j}
e^{-2\pi i\alpha n} d\alpha.
\en
It follows that
\begin{eqnarray}\label{rewrite2}
\cs_n&=&e^{n\sigma} \prod_{j=1}^n \left(
{1+e^{-j\sigma}}\right)^{a_j} \int_0^1 \prod_{j=1}^n
\left( \frac{1+e^{-j\sigma}e^{2\pi i\alpha
j}}{1+e^{-j\sigma}}\right)^{a_j} e^{-2\pi i\alpha n}
d\alpha\nonumber\\ &=& e^{n\sigma} \prod_{j=1}^n
\left( {1+e^{-j\sigma}}\right)^{a_j} \int_0^1
\phis(\alpha)\, e^{-2\pi i\alpha n} d\alpha,
\end{eqnarray}
where
$$
\phis(\alpha)=\prod_{j=1}^n \phis_j(\alpha),\quad\alpha\in\BR,
$$
and
$$
\phis_j(\alpha)=\left( \frac{1+e^{-j\sigma}e^{2\pi
i\alpha j}}{1+e^{-j\sigma}}\right)^{a_j}.
$$
If $\sigma >0,$ then the $\phis_j$ are characteristic
functions of a sequence of independent binomial random
variables $j^{-1}\Xs_j:$

$$ \BP(\Xs_k=jl)={a_k\choose
l}\left(\frac{e^{-k\sigma}}{1+e^{-k\sigma
}}\right)^l
\left(\frac{1}{1+e^{-k\sigma}}\right)^{a_k-l}, \ \
l=0,1,2,\ldots,a_k. $$ The formula (\ref{rewrite2}) could also
be derived from (145) of \cite{AT}.

The number of integer partitions of $n$ with distinct parts all of size
at least
$s$ was considered in \cite{FP}. This is the selection with
$$
a_j=\left\{
\begin{array}{l l}
0&{\rm if \ }j<s,\\
1&{\rm if \ }j\geq s.
\end{array}
\right.
$$
The identity (\ref{rewrite2}) was derived in \cite{FP} for this particular
example.

Let $\Ys=\sum_{j=1}^n \Xs_j$. We have $$ \Ms_n:=\BE\Ys =
\sum_{j=1}^n\frac{ja_je^{-j\sigma}}{1+e^{-j\sigma}}.
$$ We will assume that $\sigma=\sigmas_n$ is chosen in such
a way that $\Ms_n=n$. The fact that $\sigmas_n$ can be
chosen in such a way follows from considering that $\Ms_n$
decreases from ${1\over 2}\sum_{j=1}^nja_j$ to $0,$ as
$\sigma$ changes from $0$ to $+\infty$ and noting that the
assumption that the $a_j$ satisfy (\ref{exp}) implies that
$\sum_{j=1}^nja_j>n$ for $n$ large enough. Under the
above choice of $\sigma,$ the variance of $\Ys$ is $$
(\Bs_n)^2:={\rm Var}(\Ys)=
\sum_{j=1}^n\frac{j^2a_je^{-j\sigmas_n}}{(1+e^{-j\sigmas_n})^2}.
$$

From this starting point the
proof of Theorem \ref{selection} is similar to the proofs
of Theorem~\ref{main} and Theorem~\ref{firstorder}.


\begin{thebibliography}{99}

\bibitem{Andrews} G. Andrews,
The Theory of Partitions, Cambridge University Press, 1998.

\bibitem{ABT} R. A. Arratia, A.Barbour and S. Tavar\'e,
Logarithmic combinatorial structures: A probabilistic
approach, EMS monograph in Mathematics, vol.1, European
Mathematical Society Publishing house, Zurich, 2003.

\bibitem{AT} R. A. Arratia and S. Tavar\'e,
Independent process approximations for random combinatorial
structures, {\em Adv. Math.}, {\bf 104} (1994), 90 -- 154.
\bibitem{Ba} L. B$\acute{a}$ez-Duarte, Hardy-Ramunujan's
asymptotic formula for partitions and the central limit
theorem. {\em Adv. Math.}, {\bf 125}, (1997), 114-120.

\bibitem{BG} A. D. Barbour and B. L. Granovsky,
Random combinatorial structures: the convergent case,{\em
J. Comb. Theory, Ser. A}, 109 (2005), 203--220.

\bibitem{Bl} J. P. Bell, sufficient conditions for
zero-one laws, {\em Trans. Amer. Math. Soc.}, {\bf 354}
(2002), 613--630.

\bibitem{BB} J. P. Bell and S. N. Burris,
Asymptotics for logical limit laws: when the growth of the
components is in an RT class, {\em Trans. Amer. Math.
Soc.}, {\bf 355} (2003), 3777 -- 3794.

\bibitem{Br} N.A. Brigham, On a certain weighted partition
function, {\em Proc. Amer.Math. Soc.}, {\bf 1}(1950),
192-204.
\bibitem{Bur} S. Burris, Number theoretic density and logical limit
laws, Mathematical surveys and monographs,{\bf 86},
American Mathematical Society, Providence, RI, 2001.

\bibitem{Ca} K. J. Compton, A logical approach to asymptotic combinatorics I.
First order properties. {\em Adv. Math.} {\bf 65} (1987) 65
-- 96.

\bibitem{Cb} K. J. Compton, A logical approach to asymptotic combinatorics II.
Monadic second-order properties. {\em J. Comb. Theory, Ser.
A} {\bf 50} (1989) 110 -- 131.

\bibitem{EG} M. M. Erlihson and B. L. Granovsky,
Reversible coagulation-fragmentation processes and random
combinatorial structures: asymptotics for the number of
groups,{\em  Random sructures and algorithms}{\bf 25}
(2004), 227-245

\bibitem{FP} G. A. Freiman and J. Pitman,
Partitions into distinct large parts, {\em J. Austral.
Math. Soc.} {\bf 57} (1994) 386 -- 416.

\bibitem{FG} G. A. Freiman and B. L. Granovsky,
Asymptotic formula for a partition function of reversible
coagulation-fragmentation processes, {\em Israel J. Math.}
{\bf 130} (2002) 259 -- 279.

\bibitem{FG2} G. A. Freiman and B. L. Granovsky,
Clustering in coagulation-fragmentation processes, random
combinatorial structures and additive number systems:
Asymptotic formulae and limiting laws,  {\em Trans. Amer.
Math. Soc.}{\bf 357}(2005) 2483--2507.

\bibitem{FR} B.Fristedt, The structure of random partitions
of large integers,{\em Trans.Amer.Math. Soc.}{\bf
337}(1993) 703-735.


\bibitem{GKP} R. L. Graham, D. E. Knuth, O. Patashnik,
Concrete Mathematics, Addison-Wesley, 1989.

\bibitem{GR} I. S. Gradshteyn and E. M. Ryzhik,
Tables of integrals, series, and products, 4th. ed.,
Academic Press, 1980.

\bibitem{H} G. H. Hardy. Ramanujan, 3rd ed., Chelsea Publishing Company,
1978.

\bibitem{K}  A. I. Khinchin, Mathematical foundations of quantum
statistics, Graylock Press, Albany, N.Y., 1960.

\bibitem{KKWa} A. Knopfmacher, J. Knopfmacher, R. Warlimont,
``Factorisatio numerorum'' in arithmetical semigroups,
{Acta. Arith.} {\bf 61} (1992) 327 -- 336.

\bibitem{MK} O. Milenkovic and K.J. Compton, Probabilistic
Transforms for Combinatorial Urns Models, {\em
Comb.,Probab.and Computing}{\bf 13},(2004) 645--675.

\bibitem{M} L.Mutafchiev, Large distinct part sizes in a
random integer partition, {\em Acta Math. Hungar.} {\bf
87}(2000) 47--69.



\bibitem{P} A. G. Postnikov, Introduction to Analytic Number Theory,
Translations of Mathematical Monographs, Vol. 68, American
Mathematical Society, Providence, RI, (1987).

\bibitem{R} L. B. Richmond, A general asymptotic result for
partitions, {\em Canad. J. Math.} {\bf 27} (1975) 1083 --
1091.

\end{thebibliography}
\end{document}